# Variability in Modeling the Cyclic Scheduling of an FMC Depending on the Underlying TSP Models

Hüseyin Güden[a], Mazyar Ghadiri Nejad[b], Gergely Kovács[c], Béla Vizvári[d]

[a,d]Department of Industrial Engineering, Eastern Mediterranean University, Famagusta, TRNC, Turkey
[b]Department of Industrial Engineering, Girne American University, Kyrenia, TRNC, Turkey
[c]Edutus University, Tatabánya, Hungary

[a]huseyin.guden@emu.edu.tr
[b]mazyarghadirinejad@gau.edu.tr
[c]kovacs.gergely@edutus.hu
[d] bela.vizvari@emu.edu.tr

**Abstract**

Flexible robotic cells (FRCs) are utilized to produce standardized products at a high-speed production rate, and to set up the production floors based on the rapid operating environment changes. In an FRC, there are a number of computer numerical control (CNC) machines, an input buffer, an output buffer, and a robot. The input and the output buffers contain unprocessed and finished items, respectively, whereas the robot performs the loading/unloading activities and transports the items among the machines and buffers, getting controlled by a central computer. The system repeats a cyclic schedule in its long run and the cycle time depends on the order of the robot activities. In order to maximize the efficiency of the system order of the robot activities yielding the minimum cycle time should be determined. Aim of this research is to find a novel exact model to solve this cyclic scheduling problem (CSP). These types of problems have tight relations with the Traveling Salesmen Problem (TSP). In this study, CSP of an FRC is considered, and the adaptations of four main modeling approaches of the TSP such as Dantzig-Fulkerson-Johnson (DFJ) approach, Miller-Tucker-Zemlin (MTZ) approach, Vajda's $n$-step modeling approach, and the network flow modeling approach, are studied. Only the DFJ approach could not be adapted to the CSP and its reason is discussed. The others are adapted successfully. Furthermore, similarities and differences between the considered CSP and TSP models are scrutinized, and their performances (pros and cons) are compared together, by using several numerical cases.

## 1. Introduction

One of the significant ways to improve traditional manufacturing industries is using robots. These robots, which are called industrial robots, are mainly used for assembling, painting, welding, and transporting processes in the production systems to improve quality, productivity, and workers' safety (Barenji et al., 2014). This capability of robots speeds up the process which causes to increase production rates.

This study considers a real-life FRC with a number of parallel CNC machines, located on a line to process items. In this cell, each machine is capable of performing all required processes for producing the finished items, and there is no necessity for an item to visit more than one machine of the system. When an item is processed by any of the machines, it becomes a finished item and should be put in the output buffer. The items are handled from the input buffer to the machines and after being processed, are transported from the machines to the output buffer by the robot. The flexibility of the machines in the system allows the cell to produce parts with a wide range of features. The robot is capable of transporting parts with different ranges of size and weight. This type of FRC repeats a cycle in its long run (Mosallaeipour et al., 2018). Decreasing the cycle time in such systems means increasing the production rate, which depends on the activity orders, and arises as an optimization problem in the related industries. This cyclic scheduling problem (CSP) of such FRCs has a strong relation with TSP (Gultekin et al., 2009). A more effective version of the model is discussed in Ghadiri Nejad et al. (2018b)

Based on the research performed by Orman and Williams (2007), and Williams (2013) several different integer programing models are formulated as models of TSP.

This paper is devoted to search possible modeling approaches to solve the CSP of FRCs optimally. For this purpose, the relations, similarities and differences between the CSP and the TSP are determined. Then, four main modeling approaches of the TSP are studied for adapting the CSP. Only the Dantzig-Fulkerson-Johnson (DFJ) approach could not be adapted. The reasons are presented. All the others such as Miller-Tucker-Zemlin (MTZ) approach, Vajda's $n$-step modeling approach, and the network flow modeling approach, are modified and applied on the CSP successfully. The details of each model are scrutinized, performances of the developed models are compared by using several numerical cases, and some surprising results are presented.

In the next section, an overview on the TSP is presented and the CSP is defined. Similarities and the differences between the two problems are explained. In Section 3, main approaches for modeling the TSP are considered, the original TSP models of the corresponding approaches and their adaptations to the CSP are given. In section 4, using several numerical cases the performances of the presented models are discussed. Finally, in the last section, the study is concluded.

## 2. The relation between TSP and CSP

### 2.1. An overview on the TSP

The TSP is one of the most well-known problems in industrial engineering, operations research, and related areas. In the TSP, there is a salesman and a set of cities. The salesman starts a tour from his home city, visits each of the other cities exactly once, and turns back to the home city. The aim of problem is to determine the tour having the minimum total distance. The problem is strongly NP-hard and has many application fields. Furthermore, it is considered as the root of several other problems such as vehicle routing problems, and different exact and heuristic methods are presented for it Appelgate *et al.* (2006) and Gutin and Punnen (2002). Schrijver (2009) discussed the early history of the TSP, in details. Piehler (1960) proposed a TSP based model aiming to find the minimal total delay time in an F|no-wait|X production environment of chemical industry. Similarly, Korte et al. (1990) developed the minimal total set-up time on a single machine based on the TSP approach, and solved the greatest TSP instances in chip production, which was related to minimizing the route of a laser beam. Robotka and Vizvári (2005) studied on determining the optimal route of a robot arm by using the TSP based modeling. Kota and Jarmai (2015) discussed the importance of the TSP problem and some of its application fields like transportation, distribution, and logistics. Shavarani *et al*. (2018) applied TSP in a drone delivery system. Ghadiri Nejad and Banar (2018) also applies it in aerial transportation. Figure 1 illustrates an instance of the TSP having seven nodes, as cities, and the numbers beside the edges, as the distances among the cities. The salesman starts from city 1, visits each city once, and comes back to the first city, where the aim is finding a tour with the minimum total distance.

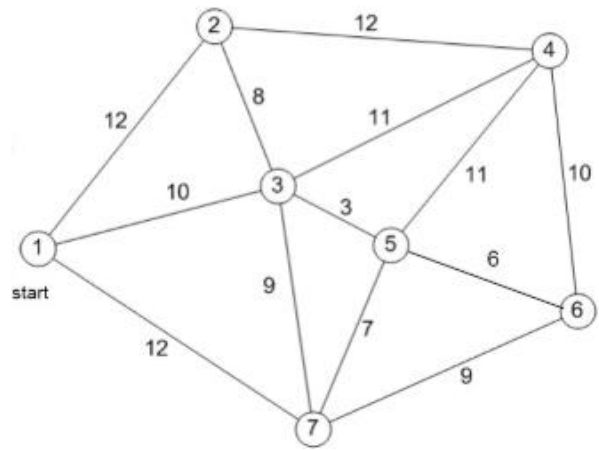

Figure 1. A TSP instance.

According to the literature, there are, approximately, ten different mathematical models for solving the TSP Orman an Williams (2007) and (Williams 2013). However, four of them are the main types of the TSP models which are significantly distinguished, such as (*i*) the DFJ model proposed by Dantzig *et al*. (1954), (*ii*) the MTZ model proposed by Miller *et al*. (1960), (*iii*) the Vajda's *n*-step model proposed by Vajda (1962), and Fox *et al.* (1980), and (*iv*) the network flow models proposed in Gavish and Graves (1978), Wong (1980) and Claus (1984).

The TSP and the famous Assignment Problem (AP) share some constraints. The main difficulty of modeling the TSP is that the optimal solution of the AP allows early return to the home city, *i.e.* the return to the home city before visiting all cities. In this way, a sub-tour is generated and the correct models of the TSP must somehow exclude the solutions consisting of two or more sub-tours. Type (*i*), *i.e.* DFJ, describes only the geometric structures of the complete tours. In general, it is the main model of the TSP. It was the first exact mathematical formulation of the TSP. It is relatively simple and is easily usable in many cases. According to Orman and Williams (2007), it is also the strongest model in a theoretical sense. Most of the numerical methods use this model Appelgate *et al*. (2006). However, it has no tool to model the time. Therefore it is not suitable for modeling scheduling problems. Each of the other three types uses a simplified description of time. However, the mathematical formulation of time is different in each case.

## 2.2. Problem definition of the considered CSP

The considered FRC consists of some parallel CNC machines located on a line. Since the machines are parallel each machine do the same manufacturing operations on the items. Hence, an unprocessed item visits only one machine and becomes a finished item. The buffers and the machines are located with the same distance from each other. Figure 2, shows a linear-layout FRC with non-identical machines. The robot moves through the line, transports the items, and performs the loading/unloading activities. If the robot goes to a machine to unload the item before the completion time of the operations of the machine it must wait there. The system repeats a cycle in its long run, in which each machine processes one part in a cycle. If the system is at a specific state at the beginning of a cycle, it reaches the same state at the end of the cycle, and then repeats the same activities in the same order in the subsequent cycles. The duration of a cycle is called cycle time. The CSP is to find the order of the loading/unloading activities which are performed by the robot in order to minimize the cycle time.

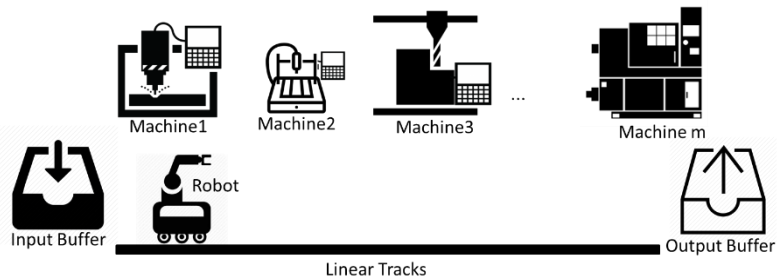

Figure 2. A robotic cell.

To formulate the problem based of different types of TSP models, the following nomenclature which are sorted alphabetically, is used: (The explanations in brackets and written italic are valid for TSP)

- $A$      Set of all loading and unloading activities performed by the robot in each cycle (*Set of the cities visited by the salesman*)
- $C$      Cycle time
- $d_{ab}$      Time of performing activity $b$ after finishing activity $a$, by the robot (*Distance from city a to city b*)
- $\delta$      Travel time of the robot between two consecutive stations
- $\varepsilon$      Time for just picking/placing an item from/to the input/output buffer, or any machine
- $L$      Set of loading activities in each cycle
- $L_i$      Loading activity of machine $i$
- $m$      Number of machines in the FRC
- $n$      (*Number of cities in the TSP instance; n=2m in the scheduling problem*)
- $p$      The processing time for an item on any machine
- $t_{ab}$      The amount of the flow from activity $a$ to activity $b$; this quantity is used only in the flow model.

| $t_a$ | Completion time of activity *a* (*The order of city a in the tour just as in MTZ model of the TSP*) |
|---|---|
| $U$ | Set of unloading activities in each cycle |
| $U_i$ | Unloading activity of machine *i* |
| $w_{ab}$ | Robot waiting time between $t_a$ and $t_b$ |
| $x_{ab}$ | 1, if the robot performs activity *b* immediately after activity *a* (*if city b is visited just after city a*); 0, otherwise |
| $v_{abs}$ | 1, if the salesman passes from city *a* to city *b* at step *s*; 0 otherwise |
| $z_i$ | 1, if $L_i$ is performed before $U_i$ for machine *i*; 0, otherwise |

The machines are assumed to be identical; thus, the process time of a part on each machine is $p$. One can use $p_i$ as process time of macgin *i* instead of $p$ if the machines are not identical. The loading activity of machine *i* ($L_i$) consists of picking, transferring, and loading a part from the input buffer to the machine *i*. Similarly, the unloading activity of machine *i* ($U_i$) includes getting the processed part from machine *i*, and transferring and putting it into the output buffer. Note that the robot stays at machine *i* at the end of any loading activity, and it stays at the output buffer at the end of any unloading activity. A cycle time is the duration spanning from the starting of the system from a specific state and returning to the same state. In order to start such a cyclic production, the system needs a setup. Each machine may be loaded or emptied at the beginning of the cycle. During a cycle, each machine must be loaded and unloaded only once. The floor plan of the FMS is as follows: The input station storing the unprocessed parts is indexed by 0. It is dollowed by the $m$ machines indexed from 1 to $m$. At the end of the line is the station of the finished pieces. It is indexed by $m + 1$. All the $m + 2$ stations and machines are in a row. The time needed by the robot to pass the distance between two neighboring station and/or robot is $\delta$. The $d_{ab}$ is calculated as follows:

$$d_{ab} = \begin{cases} 2\varepsilon + (i+j)\delta & \text{if } a = L_i \text{ and } b = L_j \\ 2\varepsilon + 2(m+1-j)\delta & \text{if } a = U_i \text{ and } b = U_j \\ 2\varepsilon + (m+1+j)\delta & \text{if } a = U_i \text{ and } b = L_j \\ 2\varepsilon + (|i-j|+m+1-j)\delta & \text{if } a = L_i \text{ and } b = U_j, i \neq j \end{cases}$$

When the process times are distinguished for the machines, some uncertain amount of waiting time for the robot can be essential. Let's consider machine *i*, its loading activity $L_i$, and unloading activity $U_i$. At the completion time of $L_i$, the machine starts its operation and finishes it after $p$ time unit. Then the robot may start unloading this part. During the unloading operation the robot takes the part from machine *i* (it takes $\varepsilon$ time units), moves to output buffer (it takes $((m + 1 - i)\delta$ time units) and put the part into the output buffer (it takes $\varepsilon$ time units). Thus, the time between the completions of $L_i$ and $U_i$ must be at least $(2\varepsilon + (m + 1 - i)\delta + p)$. There may be several other activities between $L_i$ and $U_i$, and the total time for performing those activities may not be large enough to complete the process on machine *i*. In such a case the robot must wait until the end of the process on machine *i*. The waiting times depend on the order of the activities. Note that $d_{ab}$ does not contain this uncertain amount of waiting time.

It should be noted that since the robot performs the same order of activities in a cycle, to prevent permutation and have a fix cycle, we consider $L_1$ as the last activity of the cycle. Thus, the time from finishing $L_1$ to the end of its next performing is the cycle time, and the problem is to determine the order of all loading and unloading activities in between, to minimize the cycle time.

To see an example of the sequence of robot movements in an FRC, in Figure 3 a four-machine cell has been considered, where the numbers on the arrows show the sequence of the robot movements.

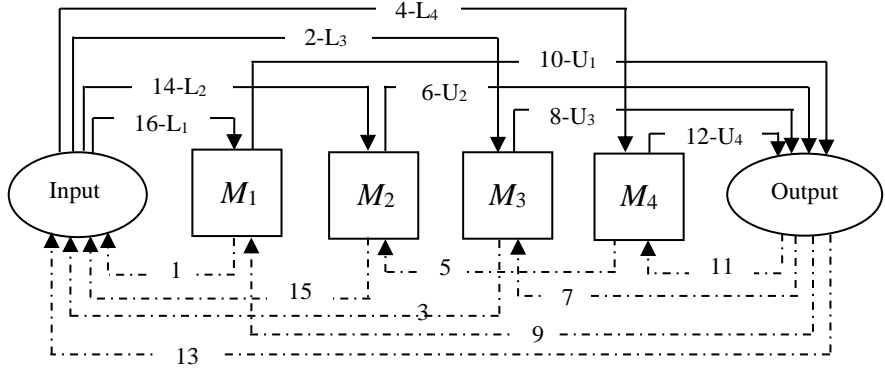

Figure 3. The sequence of robot movements for the $L_1L_3L_4U_2U_3U_1U_4L_2L_1$ cycle.

According to Figure 3 the robot is at M1 at the beginning of the cycle because the cycle starts after the completion of $L_1$. The next activity in the given order is $L_3$. Hence the robot first goes to the input buffer for taking an unprocessed part. This movement is shown by arrow 1. Then, it moves M3 to load it. This movement is shown by arrow 2.

## 2.3 The similarities and differences between the TSP and the CSP

There are some similarities between the CSP and the TSP. The robot in the CSP corresponds the salesman in the TSP; the activities in the CSP correspond the cities in the TSP; the times needed between the completions of the successive activities correspond the distances between the cities; and the cycle time in the CSP corresponds the total distance of the TSP tour. The main difference between the problems is that in the TSP, the salesman may arrive any city at any time and may leave the city at that time, while in the CSP there are pairs of activities (loading and unloading activities of the machines, $L_i$ and $U_i$ for machine $i$) that there must be at least a certain amount of time between their completions ($(2\varepsilon+(m+1-i)\delta+p)$ for machine $i$). If the robot passes the unloading activity of a pair early, then it has to wait until satisfying that certain amount of time.

## 3. TSP based models for the CSP

In order to be consistent, similar decision variables will be used in similar meaning in the TSP and the CSP models as much as possible. Let $A = \{1, 2, …, n\}$ be the set of the cities and $d_{ab}$ is the distance between city $a$ and city $b$ in the TSP network where city 1 is the home city. Since $L_1$ is considered as the first activity in the CSP it corresponds the home city of the TSP.

### 3.1. DFJ type TSP model and its adaptation to the CSP

### 3.1.1. The original DFJ type model for the TSP

The DFJ model for the TSP is the following:

$$\min \sum_{a \in A} \sum_{b \in A-\{a\}} d_{ab} x_{ab} \quad (1.1)$$

s.t

$$\sum_{a \in A-\{b\}} x_{ab} = 1 \quad \forall b \in A \quad (1.2)$$

$$\sum_{b \in A-\{a\}} x_{ab} = 1 \quad \forall a \in A \quad (1.3)$$

$$\sum_{a \in S} \sum_{b \in A-S} x_{ab} \geq 1 \quad \forall S \subset A \mid 1 \leq |S| \leq n/2 \quad (1.4)$$

$$x_{ab} \in \{0,1\} \quad \forall a,b \in A \mid a \neq b \quad (1.5)$$

In this model, objective function is the minimization of the total tour length. The first two constraints are the AP constraints and they guarantee that the salesman visits every city once. He passes from one city to another city, but without constraint (1.4) the model may generate sub-tours. Constraint (1.4) prevents sub-tours and guarantee a unique tour. By this constraint the salesman is forced to pass from a city in a subset of $A$ (which is $S$ in the constraint) to another city which is not in that subset of the cities. Since the constraint is written for all subsets of $A$, having two or more independent sub-tours is prevented. Note that there are exponential number of constraints in DFJ model.

### 3.1.2. The DFJ type model for the CSP

In the models of the CSP using only $x_{ab}$ type decision variables are not sufficient. As it is explained in the previous sections, if the robot passes from an activity to an unloading activity it can start unloading only if the part on the machine is finished. It means that if the robot arrives before the completion of the part then it must wait. The waiting time, if there is any, depends on the order of the activities. It implies that a complete description of the time is needed in any mathematical model of the problem. Hence, it also follows that DFJ is not suitable as the basis of the model.

### 3.2. MTZ type TSP model and its adaptation to the CSP

### 3.2.1. The original MTZ type model for the TSP

In the model of Miller *et al*. (1960) the home city (city 1) has a special role. According to the story of the TSP it is the city where the salesman stays. Thus, it is automatically visited. A tour can be completely described such that for every other city the order of the city in the tour is given. The order of the home city is considered 0 as no travel is needed to reach it. Hence, the orders of the other cities give the integers from 1 to $n$-1 as each city has a different order. This approach is used to eliminate sub-tours. In addition to $x_{ab}$ variables, a new type of decision variables, $t_a, \forall a \neq 1$, are used for holding the orders of the cities in the tour.

$$\min \sum_{a \in A} \sum_{b \in A - \{a\}} d_{ab} x_{ab} \qquad (2.1)$$

s.t

$$\sum_{a \in A - \{b\}} x_{ab} = 1 \qquad \forall b \in A \qquad (2.2)$$

$$\sum_{b \in A - \{a\}} x_{ab} = 1 \qquad \forall a \in A \qquad (2.3)$$

$$t_b \geq t_a + (n-1)x_{ab} - (n-2) \qquad \forall a,b \in A - \{1\} \mid a \neq b \qquad (2.4)$$

$$x_{ab} \in \{0,1\} \qquad \forall a,b \in A \mid a \neq b \qquad (2.5)$$

$$1 \leq t_a \leq n-1 \qquad \forall a \in A - \{1\} \qquad (2.6)$$

In this model sub-tours are eliminated by the constraint (2.4). The values of the order variables are restricted between 1 and $n$-1 in constraint (2.6). Let us consider two different cities $a$ and $b$. If $x_{ab}$=0 then we have $t_b \geq t_a - (n-2)$ in constraint (2.4). According to constraint (2.6) $t_a$ can be at most ($n$-1). In this case, we have $t_b \geq 1$ in constraint (2.4). Since all $t_a$ variables are greater than or equal to 1 according to constraint (2.6), if $x_{ab}$=0 constraint (2.4) becomes redundant for this ($a$, $b$) pair. On the other hand, if $x_{ab}$=1 in a solution, then the salesman passes from $a$ to $b$. The order of $b$ must be equal to the order of the city $a$ plus 1. When we substitute $x_{ab}$=1 in constraint (2.4) we obtain $t_b \geq t_a + 1$. If there is a sub-tour then there must be another sub-tour as the salesman arrives and leaves each city. Thus, if there are sub-tours then there is a sub-tour not containing the home city, *i.e.* city 1. Assume that the cities $\rho_1, \rho_2, ..., \rho_l$ are in this sub-tour, in this order. Hence,

$$x_{\rho_1 \rho_2} = x_{\rho_2 \rho_3} = \cdots = x_{\rho_{l-1} \rho_l} = x_{\rho_l \rho_1} = 1.$$

It implies that

$$t_{\rho_1} < t_{\rho_2} < \cdots < t_{\rho_l} < t_{\rho_1}$$

which is a contradiction.

Notice that the variables $t_a$ are continuous variables. The reason is that the number of these variables is $n - 1$. There are $n - 2$ differences among them. Each difference is at least 1. The values are between 1 and $n - 1$. Thus, the only possible case is that the values of the variables are the integers from 1 to $n - 1$.

The objective function and the other constraints are same with the DFJ model.

### 3.2.2. The MTZ type model for the CSP

Gültekin *et al*. (2009) published an MTZ type model. Another MTZ type model is discussed in this paper based on the model which is recently proposed by Ghadiri Nejad *et al*. (2018b). Here, the model is much more simplified and it includes some improvements. In MTZ type model of the CSP, the decision variables $x_{ab}$ have a very similar definition to the one in the original MTZ model of the TSP. The decision variables $t_{Li}$ and $t_{Ui}$, together with the new decision variables $C$ describe the cycle time. They are the equivalent of the variables of the MTZ model which describe the positions of the cities. The binary variable $z_i$ is 1 if $L_i$ precedes $U_i$ within the cycle; otherwise it is 0.

$$\min C \qquad (2.7)$$

$$\sum_{a \in A-\{b\}} x_{ab} = 1 \quad \forall b \in A) \tag{2.8}$$

$$\sum_{b \in A-\{a\}} x_{ab} = 1 \quad \forall a \in A \tag{2.9}$$

$$t_b \geq t_a + d_{ab} - M(1 - x_{ab}) \quad \forall a \neq b \in A, b \neq L_1 \tag{2.10}$$

$$t_{U_i} - t_{L_i} \leq M z_i \quad i = 1, \ldots, m \tag{2.11}$$

$$t_{U_i} \geq t_{L_i} + (2\varepsilon + (m + 1 - i)\delta + p) - M(1 - z_i) \quad i = 1, \ldots, m \tag{2.12}$$

$$t_{L_i} \leq t_{U_i} + C - (2\varepsilon + (m + 1 - i)\delta + p)(1 - z_i) \quad i = 1, \ldots, m \tag{2.13}$$

$$C \geq t_a + d_{aL_1} x_{aL_1} \quad \forall a \in A - \{L_1\} \tag{2.14}$$

$$x_{ab} \in \{0,1\} \quad \forall a \neq b \in A \tag{2.15}$$

$$t_a \geq 0 \quad \forall a \in A \tag{2.16}$$

$$C \geq 0 \tag{2.17}$$

$$z_i \in \{0,1\} \quad i = 1, \ldots, m \tag{2.18}$$

In the above model, the objective function is minimization of the cycle time. The first two constraints are equivalent to the AP constraints of the original TSP model. However, the interpretations are a bit different according to the definitions of the new problem. By these constraints, it is guaranteed that the robot performs all the activities. It passes from one activity to another activity.

The next set of constraints, constraint (2.10), is about the completion times of activities which are immediately after each other. As usual, $M$ denotes a great positive number. No exact equation is claimed in the constraints. It is not necessary for two reasons. One is that if the robot must wait at a machine, which is later in the tour, then the completion of some previous activities might be shifted in a time interval without losing the optimality. The other reason is that if equation is necessary, then it will be forced by optimality. For an activity $b$, as it is not known which activity will be the previous one, the inequality must be claimed for all potential previous activities $a$ in a way that if $a$ is not the previous one then the inequality is automatically satisfied. The so-called "Big $M$" technique is used for this purpose. If activity $b$ is not performed just after activity $a$, i.e. $x_{ab}$=0, the corresponding element of constraint (2.10) becomes $t_b \geq t_a - M$. Since the left hand side is negative for any feasible value of $t_a$, this constraint is redundant because of constraint (2.16). On the other hand, if $x_{ab}$=1 then we have $t_b \geq t_a + d_{ab}$ in constraint (2.10) which guarantees the minimum time difference ($d_{ab}$) between the completions of $a$ and $b$.

Constraints (2.11), (2.12) and (2.13) are related to the main difference between the TSP and the CSP: for any machine the time difference between the completion time of its loading activity and completion time of its unloading activity is a certain amount of time that must be spent, otherwise the robot has to wait. For machine $i$, the time difference between $t_{Li}$ and $t_{Ui}$ must be at least $(2\varepsilon + (m + 1 - i)\delta + p)$. Therefore, in the model it is essential to check the time difference between this pair of activities and it guarantees that this time difference satisfies the minimum required amount.

In a cycle, for machine $i$, if $L_i$ is earlier than $U_i$ then a part is loaded to machine $i$ and then it is unloaded in the same cycle. Thus, the time between the completions of $L_i$ and $U_i$ ($t_{Ui}$-$t_{Li}$) must be at least $(2\varepsilon + (m + 1 - i)\delta + p)$.

If $U_i$ is before $L_i$ in a cycle then it means that the loaded part is unloaded in the next cycle. In this case the total time between the load of the part and its unload, is the time from $t_{Li}$ to the end of the cycle ($C$-$t_{Li}$) plus the time from the beginning of the next cycle to $t_{Ui}$ ($t_{Ui}$-0), which is in total ($C$-$t_{Li}$+$t_{Ui}$). Hence, this time difference must be at least $(2\varepsilon + (m + 1 - i)\delta + p)$. Figure 4 illustrates the issue for $m=3$ case. In the order of Figure 4, $L_3$ is earlier than $U_3$ in a cycle. Thus, the same part is loaded and unloaded in the same cycle. On the other hand, $U_2$ is earlier than $L_2$ in a cycle, and the part loaded to machine 2 in a cycle is unloaded in the next cycle.

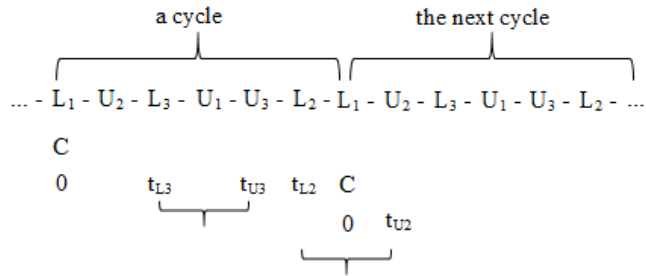

Figure 4. An example solution for $m=3$ case.

If $L_i$ is earlier than $U_i$ then $z_i$ must be 1 and ($t_{Ui}$ - $t_{Li}$) must be greater than or equal to $(2\varepsilon + (m + 1 - i)\delta + p)$. If $U_i$ is earlier than $L_i$ then $z_i$ must be 0 and ($C$-$t_{Li}$+$t_{Ui}$) must be greater than or equal to $(2\varepsilon + (m + 1 - i)\delta + p)$.

When $L_i$ is earlier than $U_i$, ($t_{Ui}$-$t_{Li}$) is positive and constraint (2.11) forces $z_i$ to be 1. Substituting $z_i$=1 in Constraints (2.12) and (2.13) yields:

$t_{Ui} \geq t_{Li} + (2\varepsilon + (m+1-i)\delta + p)$

$t_{Li} \leq t_{Ui} + C$

Therefore, constraint (2.13) becomes redundant and constraint (2.12) guarantees to have at least $(2\varepsilon + (m+1-i)\delta + p)$ time unit between $t_{Ui}$ and $t_{Li}$.

When $U_i$ is earlier than $L_i$, ($t_{Ui}$-$t_{Li}$) is negative and from the point of constraint (2.11) $z_i$ can be 0 or 1. If $z_i$=0 (which is correct according to the definition of $z_i$ and the case considered here) then constraints (2.12) and (2.13) become:

$t_{Ui} \geq t_{Li} + (2\varepsilon + (m+1-i)\delta + p) - M$

$t_{Li} \leq t_{Ui} + C - (2\varepsilon + (m+1-i)\delta + p)$

In this case, constraint (2.12) is redundant and constraint (2.13) guarantees that

$(C - t_{Li} + t_{Ui}) \geq (2\varepsilon + (m+1-i)\delta + p)$.

When $U_i$ is earlier than $L_i$, $(t_{Ui}-t_{Li})$ is negative and from the point of constraint (2.11) $z_i$ can be 0 or 1. If $z_i=1$ (which will be wrong according to the definition of $z_i$ and the case considered here) then constraints (2.12) and (2.13) become:

$$t_{Ui} \geq t_{Li} + (2\varepsilon + (m+1-i)\delta + p)$$
$$t_{Li} \leq t_{Ui} + C$$

Now, constraint (2.13) is redundant and constraint (2.12) force $t_{Ui}$ to be bigger than $t_{Li}$ which contradicts to the considered case. Therefore, when $U_i$ is earlier than $L_i$, $z_i$ can be only 0.

The cycle is finished when the robot completes $L_1$. Thus, the time from the completion of the last activity before completion of $L_1$ must be added to total time for computing $C$, correctly, which is accomplished by constraint (2.14). If activity $a$ is the mentioned activity, then, $d_{aL1}$ will be added to the cycle time. Since the objective function is the minimization of $C$ the smallest possible $C$ will be counted as the cycle time at the optimal solution. The remaining constraints are the technical constraints.

**Modification of the model to hold the waiting times**

In order to compute the waiting times, the $w_{ab}$ which is the time that the robot waits before starting activity $b$ when it is performed just after activity $a$, is considered. Then, applying the following simple modifications gives the solution with the waiting times: replace constraint (2.10) by constraints (2.19) and (2.20), replace constraint (2.14) by constraints (2.21) and (2.22), and add constraints (2.23) and (2.24) as follows:

$$t_b \geq t_a + d_{ab} + w_{ab} - M(1-x_{ab}) \quad \forall b \in A-\{L_1\}, a \in A-\{b\} \quad (2.19)$$
$$t_b \leq t_a + d_{ab} + w_{ab} + M(1-x_{ab}) \quad \forall b \in A-\{L_1\}, a \in A-\{b\} \quad (2.20)$$
$$C \geq t_a + d_{aL_1} + w_{aL_1} - M(1-x_{aL_1}) \quad \forall a \in A-\{L_1\} \quad (2.21)$$
$$C \leq t_a + d_{aL_1} + w_{aL_1} + M(1-x_{aL_1}) \quad \forall a \in A-\{L_1\} \quad (2.22)$$
$$w_{ab} \leq Mx_{ab} \quad \forall a \in A, b \in A-\{a\} \quad (2.23)$$
$$w_{ab} \geq 0 \quad \forall a \in A, b \in A-\{a\} \quad (2.24)$$

Constraints (2.19) and (2.20) compute the completion times of consecutive activities and the waiting times between them activities. Constraints (2.21) and (2.22) compute the cycle time and the waiting time between the last activity of the cycle and $L_1$. Constraint (2.23) fixes $w_{ab}$ to zero if $a$ and $b$ are not successive activities. Constraint (2.24) is the non-negativity constraint.

After the above modifications the objective function can be replaced by the following function which will give the same solution.

$$\min \sum_{a \in A} \sum_{b \in A-\{a\}} d_{ab} x_{ab} + \sum_{a \in A} \sum_{b \in A-\{a\}} w_{ab} \quad (2.25)$$

Thus, the new model is

$$(2.7) \text{ or } (2.25)$$

s.t.
(2.8), (2.9), (2.11)-(2.13), (2.15)-(2.24)

## 3.3. Vajda's *n*-step model for the TSP and its adaptation to the CSP

### 3.3.1. The original Vajda's *n*-step model for the TSP

Vajda noticed that a tour can be described by an order of the cities (1962). At the beginning, the salesman is at the home city, city 1. Each time (or step) the salesman passes from one city to another city and after the last city he passes to the home city and completes the tour. Therefore, in *n* steps the tour is completed. He used this observation to model the TSP. Vajda's model has only one type decision variable as $v_{abs}$ which is 1, if the salesman passes from city *a* to city *b* at step *s* and it is 0 otherwise.

$$\min \sum_{a \in A} \sum_{b \in A-\{a\}} d_{ab} \sum_{s=1}^{n} v_{abs} \tag{3.1}$$

s.t

$$\sum_{b \in A-\{a\}} \sum_{s=1}^{n} v_{abs} = 1 \qquad \forall a \in A \tag{3.2}$$

$$\sum_{a \in A-\{b\}} v_{abs} = \sum_{k \in A-\{b\}} v_{bk,s+1} \qquad \forall b \in A-\{1\}, s=1,\ldots,n-1 \tag{3.3}$$

$$\sum_{a \in A-\{b\}} \sum_{s=1}^{n} v_{abs} = 1 \qquad \forall b \in A \tag{3.4}$$

$$\sum_{a \in A} \sum_{b \in A-\{a\}} v_{abs} = 1 \qquad s=1,\ldots,n \tag{3.5}$$

$$v_{abs} \in \{0,1\} \qquad \forall a,b \in A \mid a \neq b, s=1,\ldots,n \tag{3.6}$$

where the objective function is the minimization of the total tour length. The first set of constraints says that the salesman leaves every city. The second set of constraints claims that if the salesman arrives to a city at step *s*, then it must leave at step *s*+1.

Following two sets of constraints are used sometimes in the model however both can be deduced from the constraints discussed above. The third set of constraints claims that the salesman arrives to every city from only one other city. The fourth set claims that it every step the salesman goes from only one city to only one another city. The technical constraint is obvious. This model is a pure 0-1 programming problem. Its drawback is that the number of variables is $n^3$. The high number of variables restricts the applicability of the model.

### 3.3.2. The Vajda's *n*-step type model for the CSP

$$\min C \tag{3.7}$$

$$\sum_{b \in A-\{a\}} \sum_{s=1}^{2m} v_{abs} = 1 \quad \forall a \in A \tag{3.8}$$

$$\sum_{a \in A-\{b\}} v_{abs} = \sum_{k \in A-\{b\}} v_{bk,s+1} \quad \forall b \in A - \{L_1\}, \quad s=1,2,\ldots,2m-1 \tag{3.9}$$

$$\sum_{a \in A-\{b\}} \sum_{s=1}^{2m} v_{abs} = 1 \quad \forall b \in A \tag{3.10}$$

$$\sum_{a \in A} \sum_{b \in A-\{a\}} v_{abs} = 1 \quad s = 1,2,\ldots,2m \tag{3.11}$$

$$\sum_{a \in A-\{L_1\}} v_{aL_1,2m} = 1 \tag{3.12}$$

$$t_b \geq t_a + d_{ab} - M\left(1 - \sum_{s=1}^{2m} v_{abs}\right) \quad \forall b \in A - \{L_1\}, \forall a \in A - \{b\} \tag{3.13}$$

$$\sum_{a \in A-\{U_j\}} \sum_{s=1}^{2m} sv_{aU_js} - \sum_{a \in A-\{L_j\}} \sum_{s=1}^{2m} sv_{aL_js} \leq (2m-1)z_j \quad j = 2,3,\ldots,2m \tag{3.14}$$

$$z_1 = 1 \tag{3.15}$$

$$t_{U_j} \geq t_{L_j} + (2\varepsilon + (m-j+1)\delta + p) - M(1-z_j) \quad j = 1,2,\ldots,m \tag{3.16}$$

$$t_{L_j} \leq t_{U_j} + C - (2\varepsilon + (m-j+1)\delta + p)(1-z_j) \quad j = 1,2,\ldots,m \tag{3.17}$$

$$C \geq t_a + d_{aL_1} v_{aL_1,2m} \quad \forall a \in A - \{L_1\} \tag{3.18}$$

$$t_a \geq 0 \quad \forall a \in A \tag{3.19}$$

$$C \geq 0 \tag{3.20}$$

$$z_j \in \{0,1\} \quad j = 1,2,\ldots,2m \tag{3.21}$$

$$v_{abs} \in \{0,1\} \quad \forall a,b \in A, a \neq b, s = 1,2,\ldots,2m \tag{3.22}$$

Objective function is the minimization of the cycle time. The first four sets of constraints have very similar meanings with the ones in the original Vajda's *n*-step model for the TSP. Here, cities are replaced by activities. Constraint (3.8) guarantees that the robot passes from each activity to another activity. Constraint (3.9) claims that if the robot performs an activity at step *s*, then it must pass from that activity to another activity at step *s*+1. Constraint (3.10) claims that the robot passes to every activity from only one other activity. Constraint (3.11) claims that at every step the robot passes from only one activity to only one other activity. Constraint (3.12) fixes to pass from an activity to the first activity, $L_1$, at the last step. Constraint (3.13) guarantees to have enough time between the successive activities for completing them. Constraint (3.14) forces $z_j$ to be 1 if $L_j$ is earlier than $U_j$. Since $L_1$ is considered as the first activity in the cycle, $U_1$ is always after $L_1$. Thus, constraint (3.15) fixes $z_1$ to 1. The functions of constraints (3.16), (3.17) and (3.18) are same with the constraints (2.12), (2.13), (2.14) of the MTZ type model of the CSP.

**Modification of the model to hold the waiting times**

It is similar to the modifications presented for the MTZ type model to update the model for holding the waiting times. It is needed to use the $w_{ab}$ variables with the same meaning. Replace constraint (3.13) by the following constraints (3.23) and (3.24), replace constraint (3.18) by the following constraints (3.25) and (3.26), and add the following constraints (3.27) and (3.28).

$$t_b \geq t_a + d_{ab} + w_{ab} - M(1 - \sum_{s=1}^{2m} v_{abs}) \quad \forall b \in A - \{L_1\}, a \in A - \{b\} \quad (3.23)$$

$$t_b \leq t_a + d_{ab} + w_{ab} + M(1 - \sum_{s=1}^{2m} v_{abs}) \quad \forall b \in A - \{L_1\}, a \in A - \{b\} \quad (3.24)$$

$$C \geq t_a + d_{aL_1} + w_{aL_1} - M(1 - v_{aL_1,2m}) \quad \forall a \in A - \{L_1\} \quad (3.25)$$

$$C \leq t_a + d_{aL_1} + w_{aL_1} + M(1 - v_{aL_1,s2m}) \quad \forall a \in A - \{L_1\} \quad (3.26)$$

$$w_{ab} \leq M \sum_{s=1}^{2m} v_{abs} \quad \forall a \in A, b \in A - \{a\} \quad (3.27)$$

$$w_{ab} \geq 0 \quad \forall a \in A, b \in A - \{a\} \quad (3.28)$$

Functions of the constraints (3.23)-(3.28) are same with the functions of the constraints (2.19)-(2.24) in the MTZ type model. After the above modifications, the objective function can be replaced by the following function which give the same solution.

$$\min \sum_{a \in A} \sum_{b \in A - \{a\}} \sum_{s=1}^{2m} d_{ab} v_{abs} + \sum_{a \in A} \sum_{b \in A - \{a\}} v_{abs} \quad (3.29)$$

$$\min \sum_{a \in A} \sum_{b \in A - \{a\}} \sum_{s=1}^{2m} d_{ab} v_{abs} + \sum_{a \in A} \sum_{b \in A - \{a\}} w_{ab} \quad (3.29)$$

Thus, the new model is

(3.7) or (3.29)
s.t.
(3.8)-(3.12), (3.14)-(3.17), (3.19)-(3.28)

### 3.4. Network flow type TSP model and its adaptation to the CSP

### 3.4.1. The original network flow type model for the TSP

Network flow type TSP model is discussed in this paper based on the model which is recently proposed by Gavish and Graves (1978). This model combines the AP as base of the model and a network flow. There is a flow on the edges of the graph. Home city is the only demand node and its demand is ($n$-1) unit. Each of the other cities has 1 unit supply. Thus, the total supply by the nodes is also ($n$-1) unit. Positive flow can exist on edges only where the salesman travels. It may be said that the salesman visits the cities, collects their supplies and bring them to the home city. Hence, outflow at city 1 is zero and inflow is ($n$-1). At each of the other cities outflow is equal to the inflow plus 1. In this model, $t_{ab}$ is the amount of the flow from city $a$ to city $b$. (i.e., the total flow amount up to visiting city $b$ when passed from city $a$).

$$\min \sum_{a \in A} \sum_{b \in A-\{a\}} d_{ab} x_{ab} \tag{4.1}$$

s.t

$$\sum_{a \in A-\{b\}} x_{ab} = 1 \quad \forall b \in A \tag{4.2}$$

$$\sum_{b \in A-\{a\}} x_{ab} = 1 \quad \forall a \in A \tag{4.3}$$

$$t_{ab} \leq (n-1) x_{ab} \quad \forall a, b \in A \mid a \neq b \tag{4.4}$$

$$\sum_{b \in A-\{a\}} t_{ab} = \sum_{k \in A-\{a\}} t_{ka} + 1 \quad \forall a \in A - \{1\} \tag{4.5}$$

$$x_{ab} \in \{0,1\} \quad \forall a, b \in A \mid a \neq b \tag{4.6}$$

$$0 \leq t_{ab} \quad \forall a, b \in A \mid a \neq b \tag{4.7}$$

The objective function and the AP constraints are the same as in the case of MTZ and the original DFJ model. The objective function is the minimization of the total tour length. The first two constraints guarantee that the salesman visits every city once. He passes from one city to another city. Constraint (4.4) allows only the edges that the salesman uses in his tour to have some positive flow. Note that the flow on any edge can be at most the total supply which is $n-1$. Constraint (4.5) is the flow conservation (or balance) equation. Outflow is inflow plus 1 for all the cities but the home city.

Assume that a sub-tour exists. It implies the existence of at least one more sub-tour. Hence, there is a sub-tour not containing the home city. Assume that the sub-tour consists of $s$ cities ($s>1$). Let $f$ be the flow value entering a certain city $p$. Notice that all cities are supply nodes in the sub-tour. Thus, the flow value increases by 1 at every city. It implies that the flow value entering city $p$ when the flow comes back is $f+s$ which is a contradiction.

### 3.4.2. The network flow type model for the CSP

Variables of type $t_{ab}$ are used in the flow model, instead of variables of type $t_a$. They represent both time and flow values. The network flow model for the CSP is the following:

$$\min \sum_{a \in A-\{L_1\}} t_{aL_1} \tag{4.8}$$

$$\sum_{a \in A-\{b\}} x_{ab} = 1 \quad \forall b \in A \tag{4.9}$$

$$\sum_{b \in A-\{a\}} x_{ab} = 1 \quad \forall a \in A \tag{4.10}$$

$$t_{ab} \leq M x_{ab} \quad \forall a \neq b \in A \tag{4.11}$$

$$w_{ab} \leq p x_{ab} \quad \forall a \neq b \in A \tag{4.12}$$

$$\sum_{b \in A-\{a\}} t_{ab} = \sum_{k \in A-\{a\}} t_{ka} + \sum_{b \in A-\{a\}} w_{ab} + \sum_{b \in A-\{a\}} d_{ab} x_{ab} \quad \forall a \in A - \{L_1\} \quad (4.13)$$

$$\sum_{b \in A-\{L_1\}} t_{L_1 b} = \sum_{b \in A-\{L_1\}} w_{L_1 b} + \sum_{b \in A-\{L_1\}} d_{L_1 b} x_{L_1 b} \quad (4.14)$$

$$\sum_{a \in A-\{U_i\}} t_{aU_i} - \sum_{a \in A-\{L_i\}} t_{aL_i} \leq M z_i \quad i = 2,3,\ldots,m \quad (4.15)$$

$$\sum_{a \in A-\{U_i\}} t_{aU_i} \geq \sum_{a \in A-\{L_i\}} t_{aL_i} + (2\varepsilon + (m-j+1)\delta + p) - M(1-z_i) \quad i = 2,3,\ldots,m \quad (4.16)$$

$$\sum_{a \in A-\{L_i\}} t_{aL_i} \leq \sum_{a \in A-\{U_i\}} t_{aU_i} + \sum_{a \in A-\{L_1\}} t_{aL_1} + (2\varepsilon + (m-j+1)\delta + p)(1-z_i) \quad i = 2,\ldots m \quad (4.17)$$

$$\sum_{a \in A-\{U_1\}} t_{aU_1} \geq p + 2\varepsilon + m\delta \quad (4.18)$$

$$t_{ab}, w_{ab} \geq 0 \quad \forall a \neq b \in A \quad (4.19)$$

$$z_i \in \{0,1\} \quad i = 2,\ldots,m \quad (4.20)$$

$$x_{ab} \in \{0,1\} \quad \forall a \neq b \in A \quad (4.21)$$

The objective is to minimize the cycle time which is the time that the robot performs $L_1$ after the last activity. Note that the cycle starts at the time that $L_1$ is completed. That time is considered as time zero. During the cycle all the activities, including $L_1$, must be completed. So, the end of a cycle is the completion time of $L_1$, which is also the beginning of the next cycle. Constraints (4.9) and (4.10) are the same constraints with the same meanings in the MTZ type model of the CSP. Constraints (4.11) and (4.12) fix $t_{ab}$ and $w_{ab}$ variables to zero if the robot does not perform activity $b$ just after activity $a$. If activity $b$ is performed just after activity $a$ then $x_{ab}$ is 1 and the corresponding $t_{ab}$ and $w_{ab}$ variables are allowed to be positive by constraints (4.11) and (4.12). Note that the waiting time for unloading a part cannot be more than the processing time. Because of this in constraint (4.12) $p$ is used as the coefficient of $x_{ab}$ instead of a big number $M$. Similar to the network flow type model of the TSP, the above constraint (4.13) is the balance constraint. Constraint (4.14) is the balance constraint for $L_1$. Constraints (4.15), (4.16) and (4.17), together, guarantee to have enough time between the load of a part and unload of it for finishing its process. Constraint (4.18) does the same thing for the part processed on the first machine.

Note that in this model the waiting times are a part of the flows between the activities and they are used in the balance constraints.

The conditions of the flow model and of the MTZ model correspond to each other according to their meaning.

## 4. Numerical results

Instead of generating new problems, in this section, all the examples considered by Ghadiri Nejad et al. (2018a) are considered, where the process times varies from zero to 250 time units, and the values of ε and δ are 1 and 2 time units, respectively. Optimization software of CPLEX 12.6 is used as the solver program and an Intel(R) Core(TM) i5-3320 CPU at 2.60GHz GHz with a RAM of 4.0 GB computer is used for the runs.

Instead of using an arbitrary very big number the cycle time of the order $L_1L_2...L_mU_1U_2...U_m$ is used as the Big $M$ value in the models. It is expected that this order gives a tight $M$ value. The cycle time of this order is computed mathematically and found as

$$M = 2(m^2 + 2m - 1)\delta + 4m\varepsilon + \max\{0, (p - 2(m-1)\varepsilon - (m^2 + m - 2)\delta)\}.$$

Tables 1 shows the optimum cycle times, solution times of the presented models (in seconds) and their linear programming relaxation (LPR) values for a set of processing time ($p$) and for $m$=4 cases.

Table 1. Solution times and LPR values for $m = 4$ cases

| $p$ | C-opt | Big $M$ Value | MTZ Type Model | | Vajda Type Model | | Network Flow Type Model | |
|---|---|---|---|---|---|---|---|---|
| | | | Solution Time (seconds) | C-LPR | Solution Time (seconds) | C-LPR | Solution Time (seconds) | C-LPR |
| 0 | 96 | 108 | 0.73 | 1.683 | 2.62 | 1.683 | 0.09 | 96 |
| 25 | 96 | 108 | 0.56 | 10.598 | 2.28 | 4.373 | 0.15 | 96 |
| 50 | 96 | 116 | 0.51 | 23.2 | 0.62 | 11.6 | 0.23 | 96 |
| 75 | 99 | 141 | 0.45 | 35.754 | 0.57 | 18.874 | 0.59 | 96 |
| 100 | 124 | 166 | 0.45 | 48.291 | 0.46 | 26.509 | 0.75 | 96 |
| 125 | 149 | 191 | 0.48 | 60.818 | 0.48 | 34.348 | 1.06 | 96 |
| 150 | 174 | 216 | 0.5 | 73.34 | 0.54 | 42.312 | 1.32 | 96 |
| 175 | 199 | 241 | 0.46 | 85.856 | 0.45 | 50.359 | 1.21 | 96 |
| 200 | 224 | 266 | 0.5 | 98.37 | 0.5 | 58.465 | 1.46 | 96 |
| 225 | 249 | 291 | 0.43 | 110.881 | 0.45 | 66.612 | 1.17 | 96 |
| 250 | 274 | 316 | 0.5 | 123.39 | 0.46 | 74.791 | 1.34 | 96 |

Figure 5 shows solution times of the presented models (in seconds) for a set of processing time ($p$) and for $m$=4 cases graphically.

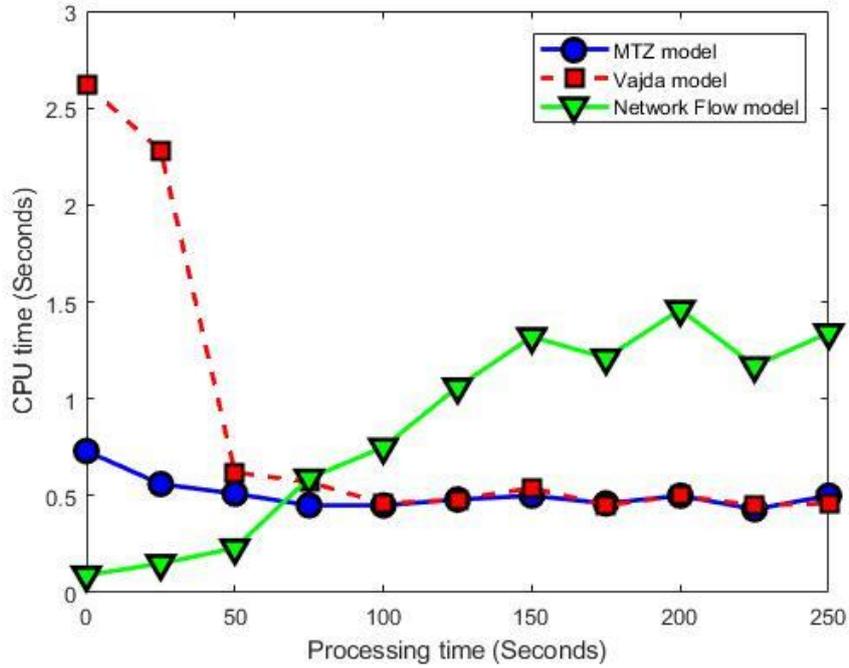
Figure 5. Solution times of the models for $m$=4 cases

According to the above results, from zero to some process time the optimal cycle time remains same. When the process time is zero it is obvious that there is no waiting time in the optimal solution. In this case, the problem reduces to a case of the TSP. The optimal solution time is equal to the sum of the robot operation times. This solution if found up to some process time. If the process time becomes great enough, then no feasible solution exists without positive waiting time. Thus, some other solutions start to be optimal. LPR of the network flow type model gives the optimal solution directly when there is no waiting time in the optimal solution. In these cases, it is the best model in terms of the solution times. Interestingly, the LPR value of the model is always same and it becomes the worst model in terms of the solution times when the optimum solutions have some positive waiting times. Then, the optimal solution has different structure. If the process time is small the LPR values of MTZ type and Vajda's $n$-step type models are worse than the flow type model but they are getting better consistently by the increase in the process time. Similarly their solution times are worse in small process times and getting better in higher process times. It seems that MTZ type model is better than Vajda's n-step type model in terms of both LPR values and solution times. Both types of models are better than the flow type model for great process times in solution times. The solution times are always less than 3 seconds, *i.e.*, very small.

Similar results are shown in Table 2, Table 3, Figure 6 and Figure 7 for $m$=5 and $m$=6 cases.

Table 2. Solution times and LPR values for $m$ = 5 cases

| $p$ | $C$-opt | Big $M$ Value | MTZ Type Model | | Vajda Type Model | | Network Flow Type Model | |
|---|---|---|---|---|---|---|---|---|
| | | | Solution Time (seconds) | $C$-LPR | Solution Time (seconds) | $C$-LPR | Solution Time (seconds) | $C$-LPR |
| 0 | 140 | 156 | 43.72 | 1.445 | 120.86 | 1.445 | 0.07 | 140 |
| 25 | 140 | 156 | 33.83 | 10.568 | 86.87 | 3.53 | 0.23 | 140 |

| | | | | | | | |
|---|---|---|---|---|---|---|---|
| **50**  | 140 | 156 | 22.4 | 23.148 | 42.85 | 9.667  | 0.21  | 140 |
| **75**  | 140 | 167 | 15.07 | 35.714 | 19.67 | 17.243 | 0.28  | 140 |
| **100** | 140 | 192 | 6.12 | 48.251 | 12.29 | 24.494 | 4.62  | 140 |
| **125** | 153 | 217 | 5.64 | 60.78  | 3.01  | 32.03  | 31.03 | 140 |
| **150** | 178 | 242 | 6.51 | 73.307 | 8.03  | 39.752 | 34.38 | 140 |
| **175** | 203 | 267 | 6.48 | 85.821 | 9.07  | 46.438 | 40.81 | 140 |
| **200** | 228 | 292 | 6.98 | 98.337 | 6.82  | 55.542 | 32.98 | 140 |
| **225** | 253 | 317 | 6.42 | 110.85 | 9.07  | 63.55  | 31.47 | 140 |
| **250** | 278 | 342 | 6.04 | 133.61 | 6.82  | 71.61  | 35.78 | 140 |

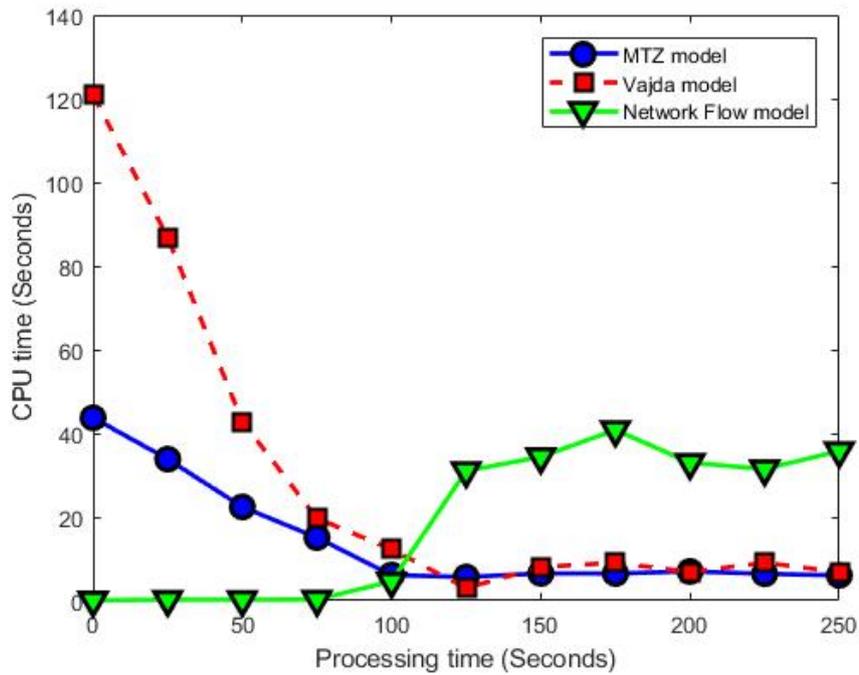

Figure 6. Solution times of the models for $m$=5 cases

Here, in $m$=5 cases, only the numbers are changed but all the explanations about the $m$=4 machine cases are valid. Solution times increased significantly. Moreover, Table 3 and Figure 7 related with $m$=6 cases support our numerical conclusions. In $m$=6 machine cases, the solution times of the Vajda's $n$-step type model is very high compared to the others. Interestingly, the LPR values of MTZ type model and Vajda's $n$-step type model are nearly the same for 4, 5 and 6 machine cases.

**Table 3.** Solution times and LPR values for $m$=6 cases

| $p$ | $C$-opt | Big $M$ Value | MTZ Type Model | | Vajda Type Model | | Network Flow Type Model | |
|---|---|---|---|---|---|---|---|---|
| | | | Solution Time (seconds) | $C$-LPR | Solution Time (seconds) | $C$-LPR | Solution Time (seconds) | $C$-LPR |
| **0** | 192 | 212 | 6330.31 | 1.289 | 37142.48 | 1.289 | 0.09 | 192 |

| 25 | 192 | 212 | 4961.29 | 10.55 | 20025.75 | 2.97 | 0.21 | 192 |
| 50 | 192 | 212 | 3686.28 | 23.109 | 7621.73 | 7.909 | 0.29 | 192 |
| 75 | 192 | 212 | 2023.67 | 35.668 | 4140.12 | 14.812 | 0.35 | 192 |
| 100 | 192 | 222 | 1402.09 | 48.217 | 2429.82 | 22.561 | 0.32 | 192 |
| 125 | 192 | 247 | 1016.54 | 60.746 | 985.52 | 29.745 | 6.48 | 192 |
| 150 | 192 | 272 | 503.28 | 73.269 | 184.78 | 37.173 | 3.71 | 192 |
| 175 | 207 | 297 | 227.18 | 85.789 | 97.66 | 44.775 | 1635.52 | 192 |
| 200 | 232 | 322 | 150.67 | 98.305 | 122.59 | 52.505 | 506.75 | 192 |
| 225 | 257 | 347 | 342.77 | 110.819 | 68.95 | 60.332 | 973.65 | 192 |
| 250 | 282 | 372 | 249.11 | 123.332 | 121.13 | 68.235 | 1231.68 | 192 |

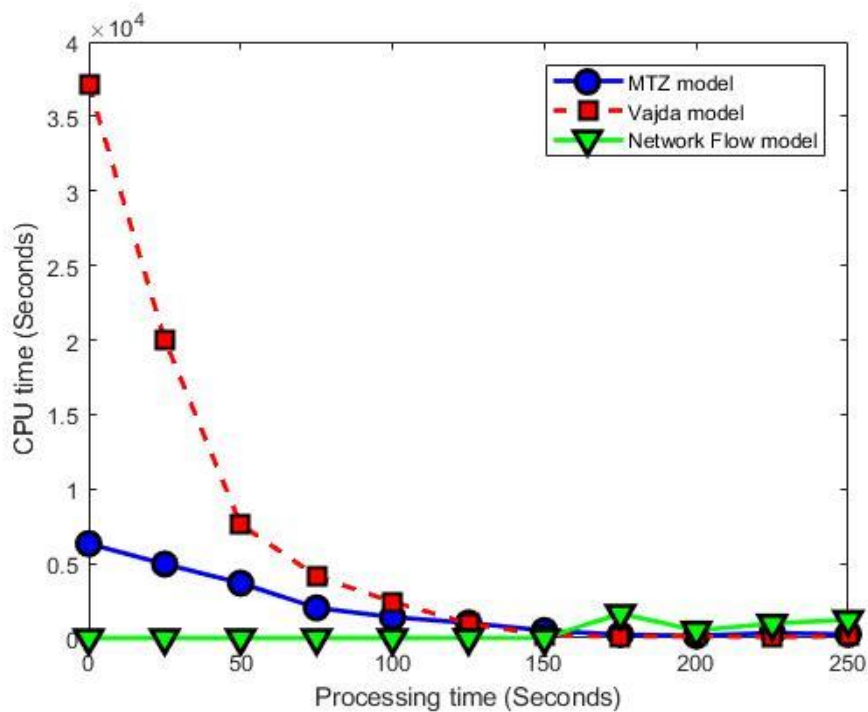

Figure 7. Solution times of the models for $m$=6 cases

## 5. Conclusion

There is a tight relation between the TSP and the CSP of flexible robotic cells. Several different types of mathematical models have been proposed to solve the TSP. In this study, the basic approaches to model the TPS have been discussed for modeling the CSP of the considered FRC. Three main approaches (MTZ, Vajda's $n$-step and Network Flow) have been successfully applied and three different models have been developed for the considered CSP. Furthermore, similarities and differences between the considered CSP and the TSP models have been scrutinized and their performances have been compared together, by using several numerical cases. The results of using linear programming relaxation of the network flow type model showed that this model finds the integer optimal solution of the problems for small process times in a very short time. At those cases, the network flow type model dominates the other methods, but when the process time is high enough to have some waiting times in the optimal solution, the results show the opposite situation.

The modelling approaches presented in this study can be used in the future to develop formulations for different versions of the CSP such as dual-gripper cases and intermediate buffer cases. New approaches may have better performances than the known approaches. Moreover, some other exact solution methods like branch-and-cut, branch-and-price can be developed based on the presented approaches. Some heuristics can be developed using these contributions.